\documentclass[letterpaper, 10 pt, conference]{ieeeconf}  % Comment this line out if you need a4paper
\IEEEoverridecommandlockouts                              % This command is only needed if 
\usepackage{xspace,empheq,fancybox,amsmath,amssymb,graphicx,epstopdf,epsfig,subfigure,syntonly,times,amsthm} 
\usepackage{psfrag,color,bm,array,cite}
\usepackage{url,cite,footnote,xspace,syntonly,algorithm,algorithmic}
\usepackage{verbatim,multirow}
\usepackage[T1]{fontenc}
\usepackage{mathtools}
\usepackage{amsmath,amsfonts,amsthm,bm} % Math packages
\usepackage{graphicx}
\usepackage{mwe}
\usepackage{caption}
\usepackage{cite}
\usepackage{amsmath,amssymb,amsfonts}
\usepackage{algorithmic}
\usepackage{graphicx}
\usepackage{textcomp}
\usepackage{xcolor}
\usepackage{multicol, blindtext}
\usepackage{xspace,empheq,fancybox,amsmath,amssymb,graphicx,epstopdf,epsfig,syntonly,times,amsthm} \usepackage{psfrag,color,bm,array,cite}
\usepackage{url,cite,footnote,xspace,syntonly,algorithm,algorithmic}
\usepackage{verbatim,multirow}
\usepackage[T1]{fontenc}
\usepackage{mathtools}
\usepackage{subfig}
\usepackage{amsmath,amsfonts,amsthm,bm}
\usepackage{graphicx}
\usepackage{mwe}
\usepackage{caption}
\usepackage{stfloats}
\usepackage{amsfonts}
\usepackage{stackengine}
\usepackage{lipsum}

\newtheorem*{remark}{Remark}

\usepackage{amsthm}
\newtheorem{proposition}{Proposition}

\newcolumntype{P}[1]{>{\centering\arraybackslash}p{#1}}
\newcommand{\hao}{\color{black}{}}
\newcommand{\yuqi}{\color{black}{}}

\newcommand*{\myov}[1]{{\overbracket[0.65pt][-1pt]{#1}}}
\renewcommand{\bar}[1]{\myov{#1}}

\title{\LARGE \bf
Substation-Level Grid Topology Optimization Using Bus Splitting}

\author{Yuqi Zhou$^{1}$, Ahmed S. Zamzam$^{2}$, Andrey Bernstein$^{2}$, and Hao Zhu$^{1}$% <-this % stops a space
\thanks{\protect\rule{0pt}{3mm} 
The work of Y. Zhou and H. Zhu was supported by NSF CAREER Grant $\#$1802319. 
This work was authored in part by the National Renewable Energy Laboratory, operated by Alliance for Sustainable Energy, LLC, for the U.S. Department of Energy (DOE) under Contract No. DE-AC36-08GO28308. Funding for A. S. Zamzam and A. Bernstein was provided by the Laboratory Directed Research and Development (LDRD) Program at NREL. The views expressed in the article do not necessarily represent the views of the DOE or the U.S. Government. The U.S. Government retains and the publisher, by accepting the article for publication, acknowledges that the U.S. Government retains a nonexclusive, paid-up, irrevocable, worldwide license to publish or reproduce the published form of this work, or allow others to do so, for U.S. Government purposes.
}
\thanks{$^{1}$Y. Zhou and H. Zhu are with the Department of Electrical \& Computer Engineering, the University of Texas at Austin, 2501 Speedway, Austin, TX 78712, USA. Emails:
        {\tt\small zhouyuqi@utexas.edu,
        haozhu@utexas.edu}}%
\thanks{$^{2}$A. Zamzam and A. Bernstein are with the Power Systems Engineering
Center, National Renewable Energy Laboratory, Golden, CO 80401, USA. Emails:
        {\tt\small ahmed.zamzam@nrel.gov, andrey.bernstein@nrel.gov}}%
}

\allowdisplaybreaks

\begin{document}

\maketitle
\thispagestyle{empty}
\pagestyle{empty}

%%%%%%%%%%%%%%%%%%%%%%%%%%%%%%%%%%%%%%%%%%%%%%%%%%%%%%%%%%%%%%%%%%%%%%%%%%%%%%%%
\begin{abstract}
Operations of substation circuit breakers are  important for maintenance needs and topology  reconfiguration in power systems. Bus splitting is one type of topology change where the two bus bars at a substation can become electrically disconnected under certain actions of circuit breakers. Because these events involve detailed substation modeling, they are typically not considered in routine power system operation and control. In this paper, an improved substation-level topology optimization framework is developed by expanding traditional line switching decisions by breaker-level bus splitting, which can further reduce grid congestion and generation costs. A tight McCormick relaxation is proposed to reformulate the bilinear terms in the resultant  optimization problem to linear inequality constraints. Thus, a tractable mixed-integer linear program reformulation is attained that allows for efficient solutions in real-time operations. Numerical studies on the IEEE 14-bus and 118-bus systems demonstrate the computational  performance and economic benefits of the proposed topology optimization approach.

%{\color{blue}{abstract not well organized, the flow can be improved}} {\hao I already edited the previous version (black), please directly modify it and add/reorganize phrases in connecting sentences}
%{\yuqi Performing topology optimization as a routine operation in the power transmission grid is becoming increasingly attractive in recent years. This technique not only can be used as corrective control actions against emergencies, but also allows for substantial operating cost reduction under the restructured system topology. 
%In this paper, a substation-level topology optimization algorithm is developed to consider both traditional line switching operations and substation bus-bar switching for enhanced system control. 
%An equivalent bus-branch model for bus splitting is first proposed through model reduction, which can be conveniently incorporated into a topology optimization framework. However, due to the potential reconnection of generators, bilinear multiplications are involved in the resultant optimization problem.
%To address this issue, we reformulate the problem using the McCormick relaxation and attain an exact mixed-integer linear program which can be efficiently solved for real-time control.
%Numerical studies on the IEEE 14-bus and 118-bus systems demonstrate the performance and considerable economic benefits of the proposed topology optimization approach.
%}

\end{abstract}


%%%%%%%%%%%%%%%%%%%%%%%%%%%%%%%%%%%%%%%%%%%%%%%%%%%%%%%%%%%%%%%%%%%%%%%%%%%%%%%%
\section{INTRODUCTION}

%{\hao (check out upper- lower- cases and line switching still present)}

Grid topology optimization is becoming increasingly important for efficient power system operations, thanks to its capability of effectively relieving network congestion and reducing generation costs.
Varying the topology of power networks mainly relies on the operations of switching devices such as circuit breakers (CBs) within the electrical substations.  %such as \textit{bus splitting}. 
The switching of CBs not only disconnects transmission lines and generation/load, but also can result in \textit{bus splitting} \cite[Ch. 11]{wood2013power}. %Operations on the CBs can give rise to system topology changes such as line switching, generation/load disconnection, and bus splitting.
%In fact, for the topology optimization problem, utilizing line switching only covers limited choices on the CB operation while more comprehensive control of breaker status (e.g., bus splitting) remains under-investigated. Therefore, in order to realize the massive potential of a comprehensive topology optimization, it is imperative to design an effective network topology optimization methodology that can further incorporate substation activities such as bus splitting.
A comprehensive topology optimization framework that includes all types of topological changes can greatly enhance the benefits  of reducing generation costs while improving the security of grid operations.  

%{\hao (first sentence issue)}
%{\yuqi The challenge of achieving substation level topology optimization is efficiently incorporating circuit breaker actions into the decision making problem.}
A majority of grid topology optimization work has mainly focused on the search of line-switching actions \cite{fisher2008optimal, hedman2009optimal, ruiz2012tractable, qiu2015chance, zhou2020transmission}, and thus have overlooked the potentials of using bus-splitting operations. %{\hao (i'm not sure what exactly is missed in the literature for the next 3 papers.)} 
The switching of substation CBs was explored in \cite{mazi1986corrective,shao2005corrective,zaoui2005coupling}  as a corrective measure for relieving localized grid stress caused by line overloads or voltage violations. %\cite{} 
%
%On the other hand, more comprehensive network topology optimization strategies that incorporate CBs activities have been proposed as corrective measures to relieve congested lines or outaged generators \cite{mazi1986corrective, shao2005corrective, zaoui2005coupling}. %account for more comprehensive system topology change, control strategies that further incorporate breaker activities are investigated thereupon. 
%In \cite{mazi1986corrective}, an algorithm was designed for corrective control by assessing the ability of candidate topology changes to relieve congestion in the network. Similarly, the corrective switching algorithm presented in \cite{shao2005corrective} incorporated more complicated bus-bar switching actions to tackle the problem of line overloads and voltage violations.
%A topology optimization model that incorporates bus splitting was presented in \cite{zaoui2005coupling} targeting network contingency analysis.
These methods have been developed to target localized contigencies in power networks by analyzing a small subset of candidate CB actions, while not yet considering a global search for the economic benefits of the full grid. 
In \cite{heidarifar2015network}, a topology optimization method based on generalized substation and CB modeling was proposed to help reduce the total generation costs. Nonetheless,  a pre-screening heuristic was utilized to address the scalability issue of the optimization problem therein to allow for real-time implementation. The optimality of the resultant topology solutions is questionable and the optimality gap is unclear. 
In fact, modeling the CB actions typically requires the detailed node-breaker representation of the power grid that includes the full list of substation components; see, e.g., \cite{pradeep2011cim, heidarifar2015network, park2019sparse,park2020optimal}. The complexity of this representation is the major cause of the lack of scalability as the resultant scenarios can be redundant. 
Thus, it still remains open to develop an efficient grid topology optimization algorithm that can account for substation-level topology change.

The goal of this paper is to develop an efficient real-time topology optimization algorithm that can incorporate the substation-level topology change such as bus splitting. %Compared with the traditional line switching operation, this can cover more comprehensive grid topology changes to relieve network congestion and reduce generation costs more effectively. 
%However, employing this model into power system operation tasks usually leads to intractable large scale optimization problems. 
To address the scalability issue of node-breaker representation, this paper leverages an equivalent bus-branch model for the substation bus splitting. Hence, instead of explicitly modeling all the components within the substation, we can conveniently incorporate this concise equivalent model for bus splitting into a grid topology optimization formulation. To deal with the bilinear terms in the resultant formulation, we apply the McCormick relaxation technique \cite{mccormick1976computability} and attain an \textit{exact} mixed-integer linear program (MILP) reformulation, which can be efficiently solved for real-time implementation.
Therefore, the main contribution of our work is to provide a tractable algorithm to effectively search for all possible topology change, both line switching and bus splitting, in order to attain the best grid-wide economic and security benefits.

The rest of the paper is organized as follows. Section \ref{sec:model} introduces the dc power flow model and the equivalent bus-branch model for bus split events.
Section \ref{sec:ots} develops the substation-level topology optimization formulation and further reformulates it into a tractable MILP using the McCormick relaxation technique.
Numerical studies on the IEEE 14-bus and 118-bus systems are presented in Section \ref{sec:num} to demonstrate the efficacy of the proposed approach. The paper is concluded in Section \ref{sec:con}.

\textit{Notation:}  Upper- (lower-) case boldface symbols are used to denote matrices (vectors); $(\cdot)^{\mathsf T}$ stands for transposition; $\mathbf I$ for identity matrix; $\mathbf 1$ denotes the all-one vector; and $\mathbf{e}_i$ denotes the standard basis vector with all entries being 0, except for the $i$-th entry which is equal to 1.

\section{System Modeling} \label{sec:model}
%The formulation of the bus split sensitivity analysis and the subsequent topology optimization formulation is based on the dc power flow model \cite{stott1974fast} and we do not consider transient stability issues for performing the network switching. 

Consider a transmission system with $N$ buses collected in the set $\cal N :=$ $\{1,\ldots,N\}$ and $L$ lines in $\cal L :=$ $\{(i,j)\} \subset \cal N \times \cal N$. For bus $i$, let $\theta_i$ be its phase angle and collect all the angles in $\bm{\theta} \in \mathbb{R}^{N}$. Similarly, let $\bm{g},~\bm{d} \in \mathbb{R}^{N}$ denote the vectors of generation and load at all buses, respectively. Under the dc power flow model, line flows $\{f_{ij}\}$ which are collected in $\bm{f} \in \mathbb{R}^{L}$ are given by: %{\hao (this sentence is a common error that I have pointed out many times. if the clause starts with a verb (utilizing), the subject of the main clause would execute this action. So you're saying lines flows are utilizing...)  }
\begin{align}
\bm{f} = \mathbf{K}\bm{\theta} \label{eq:PF1}
\end{align}
with the matrix $\mathbf{K} \in \mathbb R^{L\times N}$ mapping the phase angles to the line flows. The row of $\mathbf K$ corresponding to line $(i,j)$ is $b_{ij} (\mathbf{e}_i-\mathbf{e}_j)^\mathsf T$, where $b_{ij}$ is the inverse of the line $(i, j)$ reactance. 
%{\hao (why not in notation too?)} {\yuqi (fixed)}. 
Furthermore, the network power flow conservation leads to:
\begin{align}
\bm{p} = \mathbf{A}\bm{f} \label{eq:PF2}
\end{align}
where $\bm{p} = \bm{g} - \bm{d}$ is the net injection vector, and $\mathbf{A} \in \mathbb{Z}^{N \times L}$ is the incidence matrix for the underlying graph $(\cal N, \cal L)$.
Substituting \eqref{eq:PF1} into \eqref{eq:PF2} yields the dc power flow model:
\begin{align}
\bm{p} = \mathbf{B} \bm{\theta}
\label{eq:DCPF2}
\end{align}
where the so-termed Bbus matrix $\mathbf{B} \in \mathbb R^{N\times N}$ is given by:
\begin{align}
    \mathbf{B} = \sum_{(i,j) \in \cal L} b_{ij} (\mathbf{e}_i - \mathbf{e}_j)(\mathbf{e}_i-\mathbf{e}_j)^{\mathsf T}.
\end{align}

\begin{figure}[t!]
\centering
\vspace{2pt}
\includegraphics[trim=4.2cm 4.8cm 4cm 3cm,clip=true,totalheight=0.16\textheight]{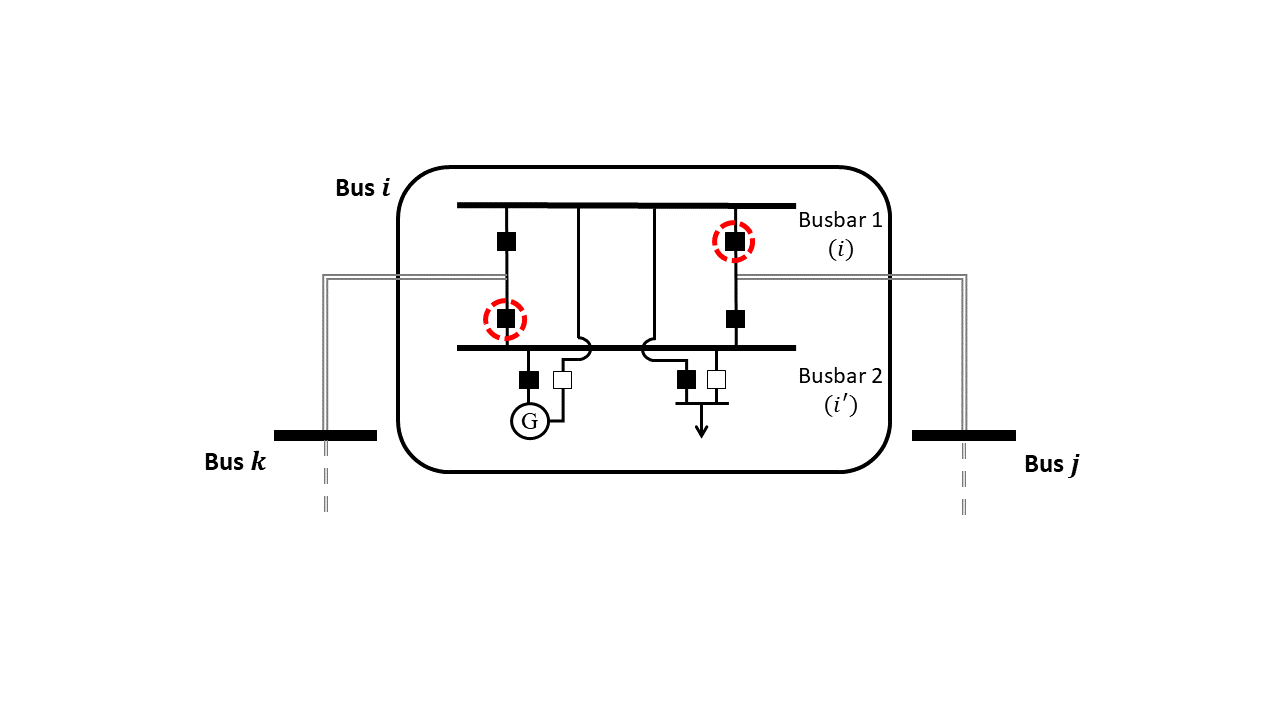}
\caption{Opening of the circled breakers leads to a bus split event at bus $i$.}
\label{fig:bus split 1}
\end{figure}

\begin{figure}[t!]
\centering
\vspace{-2pt}
\includegraphics[trim=5.5cm 7cm 7.5cm 4.8cm,clip=true,width=\linewidth]{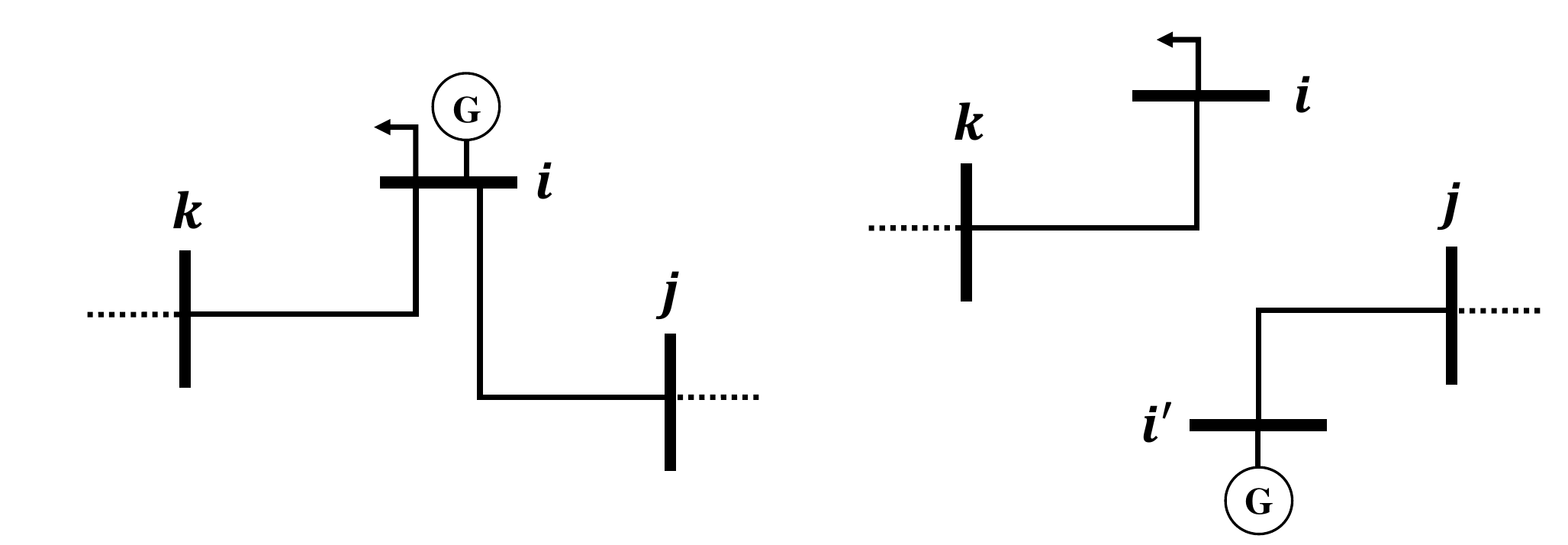}
\caption{(Left) Original bus-branch model and (right) model with a bus split at bus $i$.}
\label{fig:bus split 2}
\end{figure}

In electrical networks, switching equipment such as CBs and isolators are usually installed in substations to allow for flexible network topology and emergency intervention. 
Under certain CB configurations, the substation bus can become electrically disconnected, commonly termed as ``bus splitting'' or ``bus split.'' The occurrence of bus splitting is increasingly frequent due to misoperations of CBs \cite{kekatos2012joint, korres2006substation} or malicious cyberattacks \cite{deka2015one,ten2017impact,zhou2018false,jahromi2019cyber}.
{\yuqi Fig. \ref{fig:bus split 1} shows an example bus split event for a specific node-breaker substation configuration (double bus double breaker arrangement). Solid (hollow) squares represent closed (open) breakers. 
If the circled CBs become open, bus $i$ is split into two different buses, $i$ and $i'$. }
Accordingly, transmission lines, generation and load can be reconnected to the new bus $i'$. Although the two buses are physically co-located in the same substation, they become electrically disconnected, leading to a different bus-branch model as shown in Fig. \ref{fig:bus split 2}.

The grid-wide impact of the bus split topology change has been analyzed in \cite{zhou2019bus} and is summarized in the following proposition.

% \begin{figure}[t!]
% \centering
% \vspace{-2pt}
% \includegraphics[trim=5.5cm 6cm 9cm 5.5cm,clip=true,width=\linewidth]{bus_split_cb.pdf}
% \caption{(Left) Original substation topology and (right) the new topology with two more breakers open.}
% \label{fig:bus split 1}
% \end{figure}

\begin{proposition}
\label{prop:1}
Consider the split of bus $i$ with a single line  $(i,j)$ and injection $\tilde{p}_i$ reconnected to the new bus $i'$ (shown in Fig. \ref{fig:bus split 2}). The post-split system is equivalent to having the opening of line $(i,j)$  and an additional power transfer $\tilde{p}_i$ between buses $i$ and $j$.
%the impact on the original $N$ buses can be equivalently represented with a line outage on line $\ell=(i,j)$ and an additional power transfer $\tilde{p}_i$ between buses $i$ and $j$, where $\tilde{p}_i$ is the post-split injection at bus $i'$.
\end{proposition}

\begin{figure}[t!]
\centering
\vspace{-2pt}
\includegraphics[trim=5.5cm 6cm 7cm 6cm,clip=true,width=\linewidth]{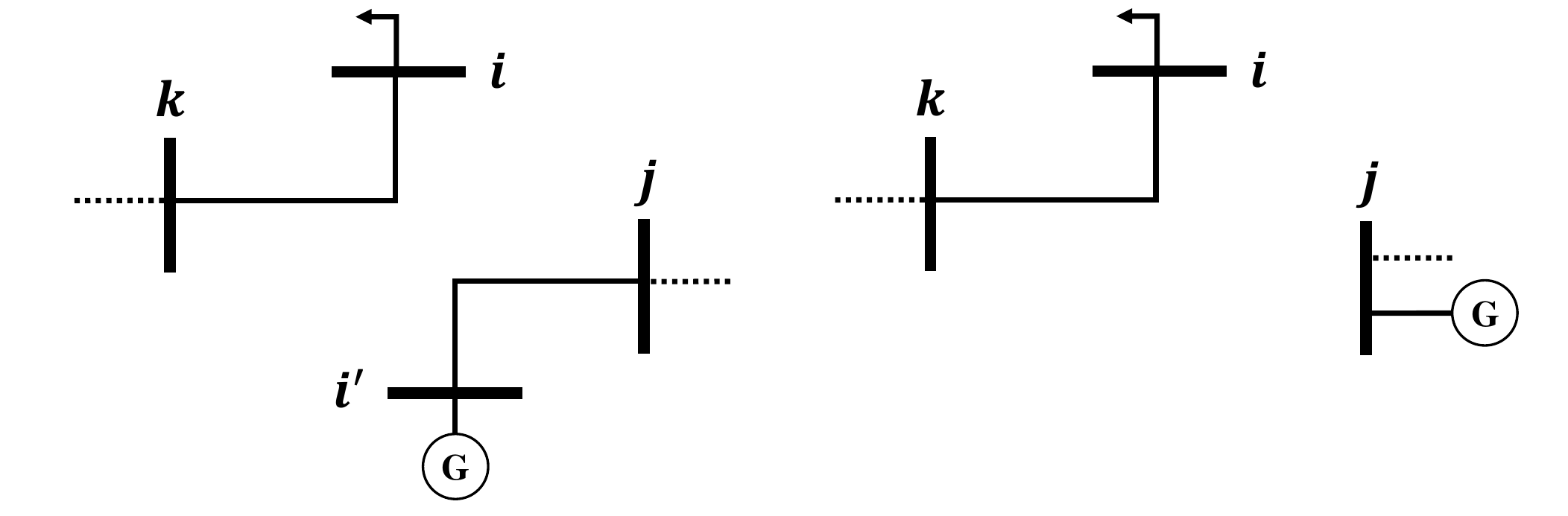}
\caption{Equivalent reduced bus-branch model (right) for the (left) system with bus split.} \label{fig:reduction}
\end{figure}

The equivalent model for the post-split system is demonstrated in Fig. \ref{fig:reduction}.
Because the new bus $i'$ connects only to bus $j$, it can be eliminated from the system by moving its connected injection ($\tilde{p}_i =g_i$ in this generation-only case) directly to bus $j$.
%, assuming no violation of transmission constraint on line $(i',j)$. 
Compared to the original system shown in Fig.\ \ref{fig:bus split 2}, the equivalent system experiences the opening of line $(i,j)$ in addition to a power transfer of $\tilde{p}_i$ between buses $i$ and $j$. This equivalent model has also been verified by linear sensitivity analysis for the bus split events \cite{zhou2019bus}. Proposition \ref{prop:1} is very useful for simplifying the incorporation of bus split events into the topology optimization problem, as discussed in the ensuing section. 

%{\hao (delete the next sentence, just want to say that 'hence' is not used as conjunction usually, not to mention favorable sounds quite strange here) Recall that traditional line switching also involves opening transmission lines, hence this equivalent model can be extremely favorable to us when being incorporated into the topology control formulation.} %The detailed formulation of substation-level topology control incorporating bus splitting is presented next.

\section{Substation Level Grid Topology Optimization} \label{sec:ots}

The grid topology optimization problem aims to determine the optimal {\yuqi transmission grid} topology with associated generation outputs in order to minimize the total generation cost. The feasible region of generation dispatch for this problem is the union of the sets of feasible solutions corresponding to each topology configuration. Thus, varying the grid topology will likely expand the overall feasible region and accordingly reduce the total generation costs \cite{fisher2008optimal, hedman2009optimal}. Going beyond the traditional topology optimization framework, the inclusion of substation level bus split events allows for additional topology change, and hence it can further reduce grid congestion and improve the security of grid operations.  

The nonlinear ac topology optimization formulation is well known to suffer from scalability issues, greatly challenging its real-time implementation \cite{soroush2013accuracies,bai2016two,zamzam2020learning}.

\begin{remark}
{\yuqi 
In this work, we adopt the dc power flow model for formulating topology optimization problem. The proposed dc based model can be possibly generalized to the nonlinear ac formulation as well; see e.g., \cite{kocuk2017new}. To corroborate the validity of the dc model, we will provide several numerical tests Section \ref{sec:num_118} to assert the feasibility of the switching solutions under the ac power flow model.
}
\end{remark}

%{\hao (delete the following but just want to point out we're repeating these sentences a lot. plz try to do more rephrasing)  However, traditional topology control focuses primarily on switching transmission lines, but fails to incorporate more complicated topology change such as bus splitting. 
%Therefore, we aim to provide the substation-level topology control which can further allow to alter the grid topology through substation bus splitting. }

\subsection{Modeling of Power Transfer in Bus Splitting}
%{\hao (the following paragraph would be a perfect case where we want the first sentence to be overarching.)}

We first discuss different scenarios of generation/load connection for modeling the power transfer in Proposition \ref{prop:1}. 
%{\yuqi Due to the flexible connectivity of generation/load after bus splitting, it is necessary to first accurately model the power transfer as stated in Proposition \ref{prop:1}.}
%The proposition states that for a single bus split, it is equivalent to a line outage and a power transfer between two buses. To conveniently incorporate the bus splitting into the topology optimization formulation, we start with the modeling of power transfer.
\begin{figure}[t!]
\centering
\vspace{6pt}
\includegraphics[trim=0cm 0cm 0cm 0cm,clip=true,totalheight=0.14\textheight]{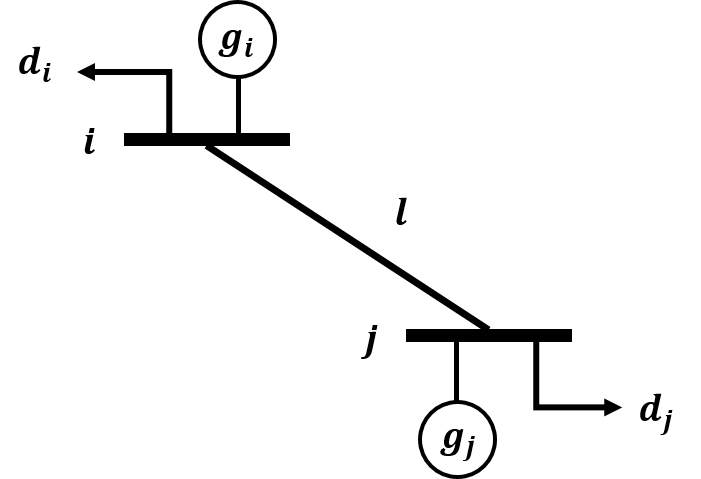}
\caption{Line $\ell = (i,j)$ with connected generation/load for modeling power transfer between buses $i$ and $j$.}\vspace{-16pt}
\label{fig:pt}
\end{figure}
To this end, consider a transmission line $\ell=(i,j)$ that connects to buses $i$ and $j$, as shown in Fig. \ref{fig:pt}. 
{\yuqi Without loss of generality,} the case of having both generation and load is assumed for the two buses. The bus split event can result in a model that is mathematically equivalent to a power transfer in between. %Essentially the power transfer involves moving one or multiple injection (generation/load) from one bus to the other one. 
For instance, the split of bus $i$ can be associated with one of three power transfer scenarios from bus $i$ to bus $j$, %that includes a total of $3$ scenarios, 
namely load only, generation only, and generation plus load. To represent the change of power injection for all three scenarios, define the following $N \times 3$  matrix:
\begin{align}
\boldsymbol{\Delta}_{\ell,i}(\bm{g}) =  (\mathbf{e}_{i} - \mathbf{e}_{j}) \left[d_i\quad -\!g_i\quad d_i\!-\!g_i\right]
\label{eq:delta_i}
\end{align}
where $\mathbf{e}_i \in \mathbb{R}^{N \times 1}$ denotes the standard basis vector.  Each column of $\boldsymbol{\Delta}_{\ell,i}(\bm{g})$ corresponds to one of three aforementioned scenarios under the split of bus $i$. Similarly, one can define this power injection matrix for the split of bus $j$, as: 
\begin{align}
{\boldsymbol{\Delta}}_{\ell,j}(\bm{g}) =  (\mathbf{e}_{j} - \mathbf{e}_{i})\left[d_j\quad -\!g_j\quad d_j\!-\!g_j\right].%\begin{bmatrix}d_j & -g_j & d_j\!-\!g_j \end{bmatrix}
\label{eq:delta_j}
\end{align}
Notice that both $\boldsymbol{\Delta}_{\ell,i}(\bm{g})$ and ${\boldsymbol{\Delta}}_{\ell,j}(\bm{g})$ depend on the generation output $\bm{g}$, which is a decision variable. In what follows, we use $\boldsymbol{\Delta}_{\ell,i}$ to refer to $\boldsymbol{\Delta}_{\ell,i}(\bm{g})$ when the dependence on $\bm{g}$ is clear from the context.

\subsection{Topology Optimization with Bus Splitting}
%{\hao (same: first sentence overarching. instead, it is repeating verbatim the first paragraph of sec III - not effective in emphasizing the key message of this subsec)}

Upon defining the power injection matrices, we are ready to formulate the topology optimization problem that includes the bus split operation.
To this end, consider the linear generation cost model, with $\bm{c} \in \mathbb{R}^{N}$ collecting all the linear coefficients. 
The binary decision variable ${z}_{\ell}$ is introduced for each transmission line $\ell=(i,j)$ to indicate the \emph{equivalent line status} (1: closed, 0: open), which will be explained in more detail after the formulation. The incident buses for line $\ell$ are collected in the set $\mathcal N_\ell$.
A vector of binary variables $\bm{w}_{\ell,i} \in \{0,1\}^{3}$ is used to select the power transfer scenario in case of a bus split at bus $i$, leading to an equivalent outage on line $\ell=(i,j)$. 
%from bus $i$ to bus $j$, for a bus split concerning line outage at line $\ell=(i,j)$.
Similarly, vector  ${\bm{w}}_{\ell,j} \in \{0,1\}^{3}$ is defined to select the power transfer scenario in case of a bus split at bus $j$, with an equivalent outage on line $\ell=(i,j)$. %from the other direction.
Under a maximum budget of $s$ operations (either line switching or bus splitting), one can formulate the following optimization problem:
\begin{subequations} \label{eq:TS1}
\begin{align}
\min \quad & {\bm{c}^{\mathsf T}\bm{g}} \label{eq:TS1_a}\\
\textrm{over} \quad  
  &  \bm{\theta} \in \mathbb{R}^{N}, \bm{g} \in \mathbb{R}^{N}, \bm{f} \in \mathbb{R}^{L}, z_{\ell} \in \{0,1\}, \; \forall \ell\in\mathcal{L} \nonumber\\
  & \bm{w}_{\ell,i} \in \{0,1\}^{3}, {\bm{w}}_{\ell,j} \in \{0,1\}^{3}, \; \forall \ell=(i,j) \in\mathcal{L} \nonumber\\
\textrm{s.t.} \quad  
  & {\theta}_{i}^{\min} \leq {\theta}_{i} \leq{\theta}_{i}^{\max},  \; \forall i \label{eq:TS1_d}\\
  & {g}_{i}^{\min} \leq {g}_{i} \leq {g}_{i}^{\max},  \; \forall i    \label{eq:TS1_e}\\
  & {f}_{\ell}^{\min} {z}_{\ell} \leq {f}_{\ell} \leq {f}_{\ell}^{\max} {z}_{\ell},  \; \forall \ell \label{eq:TS1_f}\\
  & b_{ij} (\theta_i - \theta_j) - f_{\ell} + (1-z_{\ell})\mathrm{M}_{\ell} \geq 0, \; \forall \ell = (i,j)\label{eq:TS1_g}\\
  & b_{ij} (\theta_i - \theta_j) - f_{\ell} - (1-z_{\ell})\mathrm{M}_{\ell} \leq 0, \; \forall \ell = (i,j)\label{eq:TS1_h}\\
  & \sum_{\ell}(1-z_{\ell}) \leq  s\label{eq:TS1_i}\\
  & \bm{1}^{\mathsf T} \bm{w}_{\ell,i} + \bm{1}^{\mathsf T} {\bm{w}}_{\ell,j} \leq 1 - z_{\ell},  \; \forall \ell=(i,j) \label{eq:TS1_j}\\
  & \sum_{\ell : i\in \mathcal N_\ell} \bm{1}^{\mathsf T} \bm{w}_{\ell,i}  \leq 1, \; \forall i \label{eq:TS1_k}\\
  & \mathbf{A}\bm{f}\! =\! \bm{g}\! -\! \bm{d} +\! \sum_{\ell=(i,j)}\! \boldsymbol{\Delta}_{\ell,i} \bm{w}_{\ell,i}\! +\! \sum_{\ell=(i,j)}\! {\boldsymbol{\Delta}}_{\ell,j} {\bm{w}}_{\ell,j} \label{eq:TS1_l}\\
  & f_{\ell}^{\min} \bm{1}\! \leq\!    \boldsymbol{\Delta}_{\ell,i} \bm{w}_{\ell,i} +\! {\boldsymbol{\Delta}}_{\ell,j} {\bm{w}}_{\ell,j} \leq  f_{\ell}^{\max}  \bm{1},  \; \forall \ell=(i,j) \label{eq:TS1_m}
\end{align}
\end{subequations}
%{\hao (i still believe it's easier to use sth like $\bm{w}_{\ell,i}$ $\forall i \in \mathcal N_\ell$ (i is incident to line $\ell$) instead of having the bar w. it's okay if you don't want to change it here. but definitely re-think about it when submitting the final version/preparing the slides ect.)}

%{\yuqi (The notation using $\bm{w}_{\ell,i}$ or $\bm{w}^{(i)}_{\ell}$ can cause a little confusion in constraint \eqref{eq:TS1_k}, since we are forcing that constraint for each bus $i$, and appearance of $\bm{w}^{(j)}_{\ell}$ can be difficult to understand. The concept of "direction" is critical for this constraint.)}

%{\hao (disagree. (7i) would be just $\sum_{\ell : i\in \mathcal N_\ell} \bm{1}^{\mathsf T} \bm{w}_{\ell,i}  \leq 1, \; \forall i$. Basically, it's the sum of all the w's associated with bus i, regardless of direction. )}
We discuss the constraints for problem \eqref{eq:TS1} here. Phase angle and generation limits are given in constraints \eqref{eq:TS1_d} - \eqref{eq:TS1_e}. Line flow limits are given in \eqref{eq:TS1_f}, while the flow $f_{\ell}$ is enforced to be zero when the line is open; i.e., $z_{\ell}=0$. Constraints \eqref{eq:TS1_g} - \eqref{eq:TS1_h} are introduced for establishing the line flow model in \eqref{eq:PF1}, with the constant ${\mathrm M}_{\ell}$ being sufficiently large. When the line $\ell=(i,j)$ is closed, the two inequalities are equivalent to the dc power flow equation $f_{\ell} = b_{ij} (\theta_i - \theta_j)$. Otherwise, when the line is open $f_\ell=0$ [cf. \eqref{eq:TS1_f}], the two constraints are guaranteed to be inactive for a large ${\mathrm M}_\ell$. This is called the \textit{Big-M} method \cite{griva2009linear}, which is often used to handle constraints with binary variables. For each line $\ell=(i,j)\in\mathcal{L}$, we set:
\begin{align}
\mathrm{M}_{\ell} \coloneqq b_{ij} \Delta \theta_{ij}^{\max},
\label{eq:big_M}
\end{align}
where $\Delta \theta_{ij}^{\max}$ is a given upper bound for angle stability. Constraint \eqref{eq:TS1_i} limits the total number of operations including both line switching and bus splitting, and constraint \eqref{eq:TS1_j} further defines the operations for each line.  
Specifically, if $z_{\ell}=0$ and $\bm{1}^{\mathsf T}\bm{w}_{\ell,i} + \bm{1}^{\mathsf T}{\bm{w}}_{\ell,j}=0$, the operation is simply a line switching of $\ell=(i,j)$, i.e., deenergizing the line $\ell=(i,j)$. 
Otherwise, when $z_{\ell}=0$ but $\bm{1}^{\mathsf T}\bm{w}_{\ell,i} + \bm{1}^{\mathsf T}{\bm{w}}_{\ell,j} \neq 0$, then one of the power transfer scenarios is selected after opening line $\ell$, making it equivalent to a bus split at either end of line $\ell$. {The latter case utilizes the equivalent model of bus split events and does not actually deenergize line $\ell$. Therefore,  
%are actually performing bus splitting instead of line switching, therefore 
$z_{\ell}=0$ itself cannot fully indicate the operation type (line switching or bus splitting) and is called the \emph{equivalent line status} for this reason.}

Constraints \eqref{eq:TS1_j} - \eqref{eq:TS1_m} are introduced specifically for considering bus split events.  
%Constraints \eqref{eq:TS1_j}-\eqref{eq:TS1_m} are presented specifically due to the bus split operational constraints. 
Constraint \eqref{eq:TS1_j} limits the number of power transfers that can be selected for a bus split involving the opening of line $\ell=(i,j)$. When the line is closed ($z_{\ell} = 1$), no power transfer is allowed; Otherwise, when the line is open ($z_{\ell} = 0$), at most one ($0$ or $1$) power transfer can be made, depending on whether it is line switching or bus splitting. Furthermore, notice that a single bus can be connected to multiple buses, but the power transfer from bus $i$ to other buses can be made only if the bus is split into two bus bars. Once bus $i$ is split for a power transfer with one of its incident buses, no other power transfer can be made with other buses. Therefore, constraint \eqref{eq:TS1_k} limits the power transfer from each bus to be at most once due to the physical limit of the substation. Constraint \eqref{eq:TS1_l} enforces the network power balance in \eqref{eq:PF2}, where the injection also reflects any power transfer made because of the bus split. Last, \eqref{eq:TS1_m} guarantees that for the injection reconnected to the new bus $i'$, the power flow on that incident line is not violating the transmission limit of line $\ell=(i,j)$ [cf. Fig.\ \ref{fig:reduction}].

The main challenge of solving \eqref{eq:TS1} lies in the nonlinearity of the constraints.
%The optimization problem \eqref{eq:TS1} is not in a mixed-integer linear program (MILP) formulation {\hao (this characterization strange: why would we require MILP and what is a non-tractable MILP?) } 
Specifically, constraints \eqref{eq:TS1_l} and \eqref{eq:TS1_m} are bilinear in the decision variables, due to the multiplication terms, namely $\boldsymbol{\Delta}_{\ell,i} (\bm{g})\bm{w}_{\ell,i}$ and ${\boldsymbol{\Delta}}_{\ell,j} (\bm{g}){\bm{w}}_{\ell,j}$.
To address these terms, we propose to adopt the \emph{McCormick relaxation} technique \cite{gupte2013solving} to reformulate the problem that is amenable to off-the-shelf MILP solvers.
First, rewrite the following multiplication as:
\begin{align}
\boldsymbol{\Delta}_{\ell,i} \bm{w}_{\ell,i} %& =  (\mathbf{e}_{i} - \mathbf{e}_{j}) \left[d_i(\basise_{1}+\basise_{3})^{\mathsf T} -  g_i(\basise_{2}+\basise_{3})^{\mathsf T} \right]\bm{w}_{\ell,i} \nonumber\\
& = \boldsymbol\delta_{d}^{i} \bm{w}_{\ell,i} d_i - \boldsymbol\delta_{g}^{i} \bm{w}_{\ell,i} g_i
\label{eq:multiplication_i}
\end{align}
where vectors $\boldsymbol\delta_{d}^{i} := (\mathbf{e}_{i} - \mathbf{e}_{j}) \begin{bmatrix}1 & 0 & 1\end{bmatrix}$ and $\boldsymbol\delta_{g}^{i} := (\mathbf{e}_{i} - \mathbf{e}_{j}) \begin{bmatrix}0 & 1 & 1\end{bmatrix}$.
We define the product of $\bm{w}_{\ell,i}$ and $g_{i}$ as:
\begin{align}
\bm{y}_{\ell,i} =  \bm{w}_{\ell,i} g_{i},\quad \forall~\ell. 
\label{eq:equality_1}
\end{align}
Under the bounds $[g_{i}^{\min},g_{i}^{\max}]$ for generation $g_i$, the following  four \textit{linear} inequalities hold:
\begin{subequations} \label{eq:mci}
\begin{align}
& \bm{y}_{\ell,i} \geq  \bm{w}_{\ell,i} g_{i}^{\min}\label{eq:mci_c}\\
& \bm{y}_{\ell,i} \geq \bm{1} g_{i} +  \bm{w}_{\ell,i} g_{i}^{\max} - \bm{1} g_{i}^{\max}\label{eq:mci_d}\\
& \bm{y}_{\ell,i} \leq  \bm{w}_{\ell,i} g_{i}^{\max}\label{eq:mci_e}\\
& \bm{y}_{\ell,i} \leq \bm{1} g_{i} +  \bm{w}_{\ell,i} g_{i}^{\min} - \bm{1} g_{i}^{\min}. \label{eq:mci_f}
\end{align}
\end{subequations}
The inequalities \eqref{eq:mci_c} - \eqref{eq:mci_f} can be verified by substituting \eqref{eq:equality_1}. Conversely, for any binary $\bm{w}_{\ell,i}$, the linear inequalities in \eqref{eq:mci} also guarantee the validity of \eqref{eq:equality_1}. When the $k$-th entry in the binary $\bm{w}_{\ell,i}$ is equal to zero, the two inequalities \eqref{eq:mci_c} and \eqref{eq:mci_e} jointly force the $k$-th entry of $\bm{y}_{\ell,i} $ to be {zero}. Otherwise, when the $k$-th entry of $\bm{w}_{\ell, i}$ is equal to one, the inequalities \eqref{eq:mci_d} and \eqref{eq:mci_f} enforce that the $k$-th entry of $\bm{y}_{\ell,i}$ is equal to $g_{i}$.
Due to the binary vector $\bm{w}_{\ell,i}$, %{\hao (nature???)} %{\yuqi (this word was first used in Vassilis's ACC paper; another definition of nature: the basic or inherent features of something, especially when seen as characteristic of it.)}, 
the set of inequalities in \eqref{eq:mci} is equivalent to the bilinear relation in \eqref{eq:equality_1}.
Reformulating \eqref{eq:equality_1} with the linear inequalities in \eqref{eq:mci} is known as the \textit{McCormick relaxation} technique, which has been popularly used in other problems of designing grid topology  \cite{bhela2019designing,bazrafshan2019optimal}. %{\hao ([23] used it too?)} {\yuqi (Yes)}. 

Hence, the bilinear product as given in \eqref{eq:multiplication_i} can be \textit{equivalently} replaced with:
\begin{align}
\boldsymbol{\Delta}_{\ell,i} \bm{w}_{\ell,i} = \boldsymbol\delta_{d}^{i} \bm{w}_{\ell,i} d_i - \boldsymbol\delta_{g}^{i} \bm{y}_{\ell,i}
\label{eq:multiplication_i_new}
\end{align}
and the linear inequality constraints \eqref{eq:mci_c} - \eqref{eq:mci_f}. Similarly, the bilinear product of ${\boldsymbol{\Delta}}_{\ell,j}$ and ${\bm{w}}_{\ell,j}$ can be directly given as: 
\begin{align}
{\boldsymbol{\Delta}}_{\ell,j} {\bm{w}}_{\ell,j} = \boldsymbol\delta_{d}^{j} {\bm{w}}_{\ell,j} d_j - \boldsymbol\delta_{g}^{j} {\bm{y}}_{\ell,j}
\label{eq:multiplication_j_new}
\end{align}
for similarly defined $\boldsymbol\delta_{d}^{j}$
%= (\mathbf{e}_{j} - \mathbf{e}_{i}) (\basise_{1}+\basise_{3})^{\mathsf T}$ 
and $\boldsymbol\delta_{g}^{j}$,  %= (\mathbf{e}_{j} - \mathbf{e}_{i}) (\basise_{2}+\basise_{3})^{\mathsf T}$ 
together with the following four linear inequalities:
\begin{subequations} \label{eq:mcj}
\begin{align}
& {\bm{y}}_{\ell,j} \geq  {\bm{w}}_{\ell,j} g_{j}^{\min}\label{eq:mcj_c}\\
& {\bm{y}}_{\ell,j} \geq \bm{1} g_{j} +  {\bm{w}}_{\ell,j} g_{j}^{\max} - \bm{1} g_{j}^{\max}\label{eq:mcj_d}\\
& {\bm{y}}_{\ell,j} \leq  {\bm{w}}_{\ell,j} g_{j}^{\max}\label{eq:mcj_e}\\
& {\bm{y}}_{\ell,j} \leq \bm{1} g_{j} +  {\bm{w}}_{\ell,j} g_{j}^{\min} - \bm{1} g_{j}^{\min}. \label{eq:mcj_f}
\end{align}
\end{subequations}
Thus, we have reformulated the bilinear products $\boldsymbol{\Delta}_{\ell,i} \bm{w}_{\ell,i}$ and ${\boldsymbol{\Delta}}_{\ell,j} {\bm{w}}_{\ell,j}$ using linear constraints \eqref{eq:multiplication_i_new} and \eqref{eq:multiplication_j_new} followed by additional linear inequalities \eqref{eq:mci_c} - \eqref{eq:mci_f} and \eqref{eq:mcj_c} - \eqref{eq:mcj_f}.
{\yuqi Therefore, the equivalent topology optimization problem that incorporates substation bus splitting operation can be established in the following proposition.}
\begin{proposition}
\label{prop:1}
The original nonlinear optimization problem \eqref{eq:TS1} is equivalent to the following one:
\begin{subequations} \label{eq:TS2}
\begin{align}
\min \quad & {\bm{c}^{\mathsf T}\bm{g}} \label{eq:TS2_a}\\
\textrm{over} \quad  
  &  \bm{\theta} \in \mathbb{R}^{N}, \bm{g} \in \mathbb{R}^{N}, \bm{f} \in \mathbb{R}^{L}, z_{\ell} \in \{0,1\}, \; \forall \ell\in\mathcal{L} \nonumber\\ 
  & \bm{w}_{\ell,i} \in \{0,1\}^{3}, {\bm{w}}_{\ell,j} \in \{0,1\}^{3}, \; \forall \ell=(i,j)\in\mathcal{L} \nonumber\\
  &\bm{y}_{\ell,i} \in \mathbb{R}^{3}, {\bm{y}}_{\ell,j} \in \mathbb{R}^3 \qquad \qquad\  \; \forall \ell=(i, j) \in\mathcal{L} \nonumber\\
  \textrm{s.t.} \quad  
  &   \eqref{eq:TS1_d} - \eqref{eq:TS1_k}, \eqref{eq:mci_c} - \eqref{eq:mci_f}, \eqref{eq:mcj_c} - \eqref{eq:mcj_f} \\%\label{eq:TS2_c}\\
  \begin{split}
      &\mathbf{A}\bm{f} = \bm{g} - \bm{d} + \sum_{\ell=(i, j)} (\boldsymbol\delta_{d}^{i} \bm{w}_{\ell,i} d_i - \boldsymbol\delta_{g}^{i} \bm{y}_{\ell,i})\\
      &\quad\qquad+ \sum_{\ell=(i, j)} (\boldsymbol\delta_{d}^{j} {\bm{w}}_{\ell,j} d_j - \boldsymbol\delta_{g}^{j} {\bm{y}}_{\ell,j}) \label{eq:TS2_b}
  \end{split}\\
  \begin{split}
      f_{\ell}^{\min} \bm{1} &\leq     \boldsymbol\delta_{d}^{i} \bm{w}_{\ell,i} d_i - \boldsymbol\delta_{g}^{i} \bm{y}_{\ell,i}   +   \boldsymbol\delta_{d}^{j}{\bm{w}}_{\ell,j} d_j  \\
      &\qquad \qquad\quad   - \boldsymbol\delta_{g}^{j} {\bm{y}}_{\ell,j}   \leq f_{\ell}^{\max} \bm{1}, \quad  \forall \ell=(i, j). \label{eq:TS2_c}
  \end{split} %\\
%   &  \eqref{eq:mci_c} - \eqref{eq:mci_f}, \eqref{eq:mcj_c} - \eqref{eq:mcj_f} \nonumber%\label{eq:TS2_e}
\end{align}
\end{subequations}
\end{proposition}
%Thus, we have converted the original problem \eqref{eq:TS1} into an MILP with linear constraints 
This reformulated problem is a mixed-integer linear program (MILP) and can be efficiently solved by common optimization solvers such as CPLEX, MOSEK and Gurobi.
%This \textit{equivalent} reformulation \eqref{eq:TS2} provides a tractable solution for performing the substation-level topology optimization incorporating bus split events.  
{\yuqi %In essence, the reformulation above is a tight relaxation of the original mixed-integer non-convex problem \eqref{eq:TS1}. Given that this relaxation involves binary variables and the bounds for the continuous variables are specified, the relaxed optimization problem \eqref{eq:TS2} is an \textit{exact} tight reformulation.
Moreover, since this work directly considers a minimal set of decisions on  line switching and power transfer as a result of substation bus splits, the dimension of binary decision variables has been greatly reduced from  the original one under the detailed node-breaker representation.} {\hao This reduction does not affect the  solution quality, as one can easily recover the underlying CB status and thus determine the corresponding breaker actions. Therefore, our proposed approach has utilized the concise bus-branch representation to attain an efficient and effective transmission grid switching solution using substation-level actions.  }

%{\yuqi Without compromising the solution quality, one can easily interpret and map the solutions to the underlying CB status for switching operations. Therefore, our approach \textbf{built}  upon the bus-branch representation provides an efficient transmission grid switching solution that can incorporate substation-level topology change.
%}

%%%%%%%%%%%%%%%%%%%%%%%%%%%%%%%%%%%%%%%%%%%%%%%%%%%%%%%%%%%%%%%%%%%%%%%%%%%%%%%%
\section{Numerical Results} \label{sec:num}
In this section, we first use the IEEE 14-bus system to mainly illustrate that bus split operations can be used to effectively relieve network congestion and to help address feasibility issues of the optimal power flow problem. After that, we perform the substation-level topology optimization on the larger sized IEEE 118-bus system to demonstrate the economic improvement on generation dispatch and to assess the computational complexity of the proposed topology optimization model.
The optimization problems have been implemented on a regular laptop with Intel CPU @ 2.60 GHz and 12 GB of RAM in the MATLAB R2018a simulator. %with power system modeling provided by MATPOWER \cite{zimmerman2010matpower}. 
The MILP-based optimization problems are computed using the CPLEX solver.

\subsection{14-Bus System Test}
\label{sec:num_14}
The IEEE 14-bus system consists of 20 transmission lines and 5 conventional generators, and we use the ac power flow model to test the system. The system has been slightly modified to illustrate that the bus splitting can be used to relieve network congestion and thus can help with the feasibility issue of the optimal power flow problem. Specifically, we modify the transmission limits of lines $(2, 3)$ and $(3, 4)$ to be $100$ MW and $10$ MW, respectively. The maximum generation limit of the generator at Bus 3 is adjusted to be $20$ MW, and all other network configurations are kept unaltered.

\begin{figure}[t!]
\centering
\vspace{-2pt}
\includegraphics[trim=8cm 3cm 3cm 2cm,clip=true,totalheight=0.25\textheight]{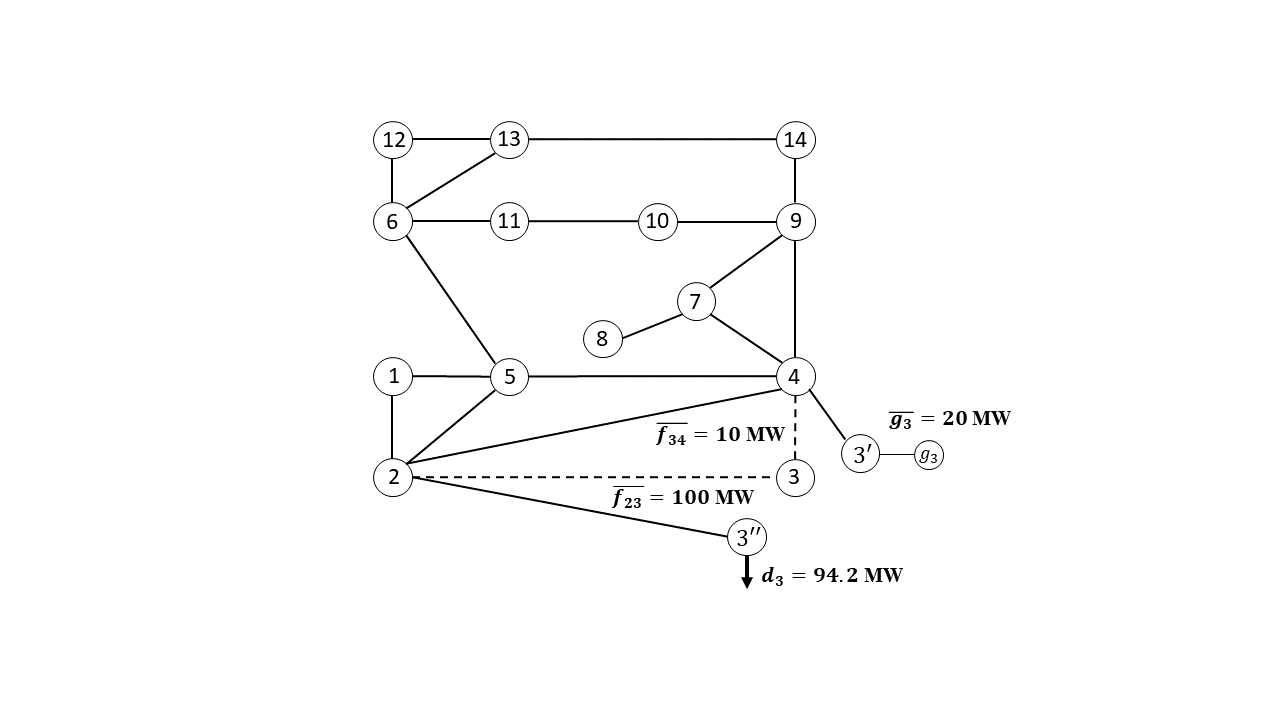}
\caption{Bus splitting at Bus 3 in the IEEE 14-bus system.}\vspace{-10pt}
\label{fig:feasibility}
\end{figure}

For Bus 3 in the original system, as shown in Fig. \ref{fig:feasibility}, a net load of at least $\bm{d}_3 - \overline{\bm{g}_3} = 74.2$ MW needs to be satisfied by the flows from line $(2, 3)$ and line $(3, 4)$. Due to the electrical characteristics of the lines, power flows on both lines are governed by the phase angle at Bus 3. Therefore, as we gradually increase the flow on line $(2, 3)$, the transmission limit on line $(3, 4)$ will be eventually violated before the sum of power flows on both lines meets the net load at Bus 3, leading to an \textit{infeasible} solution to the optimal power flow problem. 
In order to relieve the congestion on line $(3, 4)$ in this scenario, solving the proposed topology optimization problem suggests performing a bus split at Bus 3 such that the generator is connected to bus bar $3^{'}$ and the load is connected to bus bar $3^{''}$. Essentially, this bus splitting decouples the power flows on lines $(2, 3)$ and $(3, 4)$ by allowing each bus bar to have a different phase angle. 
After the bus splitting, the load at Bus 3 can be easily satisfied by the flow from line $(2, 3)$ without violating any transmission limit. Thus, the ac power flow model of the system for the updated system as shown in Fig. \ref{fig:feasibility} becomes \textit{feasible}. In fact, the bus split operation can be easily achieved through switching associated CBs at the substation. The case study on this small system indicates that similar to traditional line switching and load shedding, the operation of bus splitting can be also used to relieve network congestion and help with feasibility issues. Although the model used in \eqref{eq:TS2} is a dc power flow model, the solution obtained which suggests a bus splitting at Bus 3 makes the problem feasible while considering the ac power flow model.
Next, we will use a larger sized system to illustrate the enhanced economic benefit of the proposed substation-level topology optimization algorithm.

\subsection{118-Bus System Test}
\label{sec:num_118}
The IEEE 118-bus test case is tested for the substation-level topology optimization. The system consists of 118 buses, 186 transmission lines and 19 committed conventional generators. To illustrate the improvement and economic benefits of the topology optimization by incorporating bus split events, we also tested the same system for the traditional topology optimization strategy \cite{fisher2008optimal}. This can be easily fulfilled by restricting $\bm{w}_{\ell,i}$ and ${\bm{w}}_{\ell,j}$ in \eqref{eq:TS1} to be zero, which will exclude bus splits from consideration and allow for only line-switching operations. 
By doing so, constraints \eqref{eq:TS1_j}, \eqref{eq:TS1_k} and \eqref{eq:TS1_m} always hold and are thus disabled. Meanwhile, \eqref{eq:TS1_l} becomes a linear constraint, which describes the network power balance without any power transfer. 
Accordingly, the optimization problem \eqref{eq:TS1} itself constitutes an MILP that is readily solvable for common optimization solvers.

\begin{figure}[t!]
\centering
\vspace{2pt}
\includegraphics[trim=0cm 0cm 1cm 0cm,clip=true,totalheight=0.21\textheight]{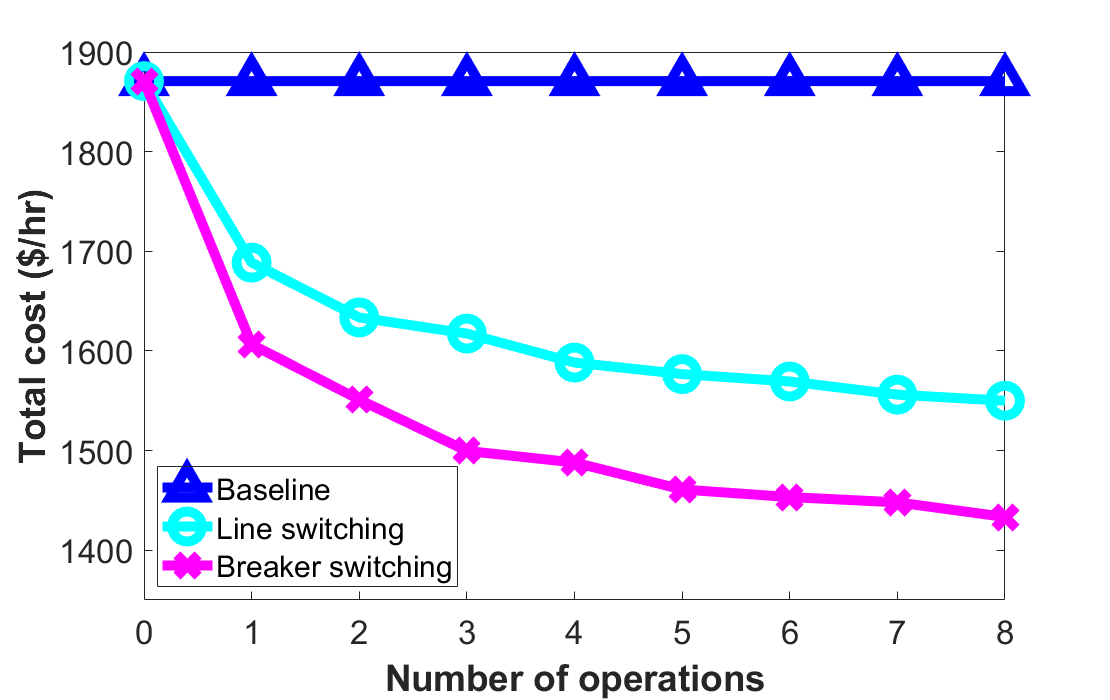}
\caption{Comparison of total cost between line switching and breaker switching for the 118-bus system.}
\label{fig:cost}
\end{figure}

The comparison of total cost under different numbers of operations for line switching and breaker-level switching is given in Fig. \ref{fig:cost} together with the benchmark cost for the system without any topology switching. Compared with the benchmark cost which involves no topology optimization, our proposed breaker-switching strategy achieves total savings of $14.1\% - 23.4\%$, depending on the number of operations; cf. Fig. \ref{fig:cost}. Meanwhile, compared with line switching, it provides additional cost savings of $4.9\% - 7.5\%$ correspondingly. Notice that these additional savings are obtained only by altering the status of several breakers at the substations, therefore the economic benefits are indeed attractive for system operators.

\begin{table}[t!] 
\caption{Comparison of Topology Optimization Decisions for the 118-Bus System} \label{tab:sd} 
\centering 
{\renewcommand{\arraystretch}{1.2}
\begin{tabular}{P{0.4cm} p{2.6cm} p{2.6cm} p{1.2cm} }
\hline 
 $s$ &  Line Switching   &  Breaker Switching  &  Reduction \\ 
\hline 
1 & Line 128 & Bus 82 & 4.9\%\\   
2 & Lines 128, 136  & Buses 77, 82 & 5.1\%\\  
3 & Lines 41, 128, 136  & Buses 77, 82 \& Line 130 & 7.3\%\\  
4 & Lines 119, 123, 124, 125  & Buses 75, 77, 82 \& Line 136 & 6.3\%\\  
5 & Lines 118, 121, 131, 135, 149  & Buses 77, 82 \& Lines 123, 124, 125 & 7.3\%\\  
\hline  
\end{tabular}}
\end{table}

\begin{figure}[t!]
\centering
\vspace{-2pt}
\includegraphics[trim=0cm 0cm 1cm 0cm,clip=true,totalheight=0.21\textheight]{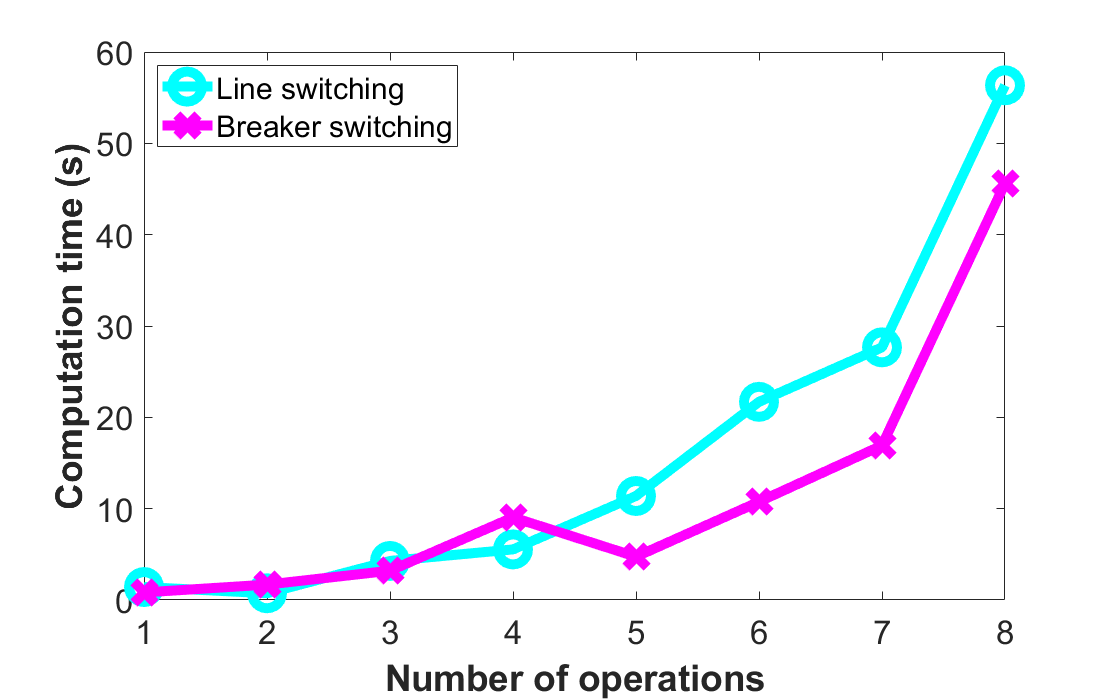}
\caption{Comparison of computation time between line switching and breaker switching for the 118-bus system.}\vspace{-14pt}
\label{fig:time}
\end{figure}

To compare the {\yuqi switching decisions} provided by the traditional line switching and the proposed breaker switching, we list the topology optimization solutions for up to a maximum of $s=5$ operations in Table \ref{tab:sd}. In the line-switching part, only operation on the transmission lines is allowed, which normally involves opening a pair of breakers at both ends of the line. In contrast, the breaker-switching strategy enables more complicated breaker operations that lead to not only line switching but also bus split events. Take $s=3$ as an example, the line-switching scheme picks lines 41, 128 and 136 to open, whereas the proposed breaker switching suggests opening line 130 and performing bus splits at buses 77, 82 simultaneously. As the result of considering breaker-level operations, an additional reduction of $7.3\%$ in the operational cost is achieved. 

Additionally, we assessed the computation time for the two formulations of the topology optimization scheme shown in Fig. \ref{fig:time}. For up to $s=8$ operations, on average the solving time for the substation-level topology optimization is $28.1\%$ faster than the line-switching one. Notice that compared with the line-switching formulation, we further introduce variables such as $\bm{w}_{\ell,i}$, ${\bm{w}}_{\ell,j}$, ${\bm{y}}_{\ell}$ and ${\bm{y}}_{\ell,j}$, but the additional constraints \eqref{eq:TS1_j} - \eqref{eq:TS1_k} that incorporate bus split events can potentially facilitate the computation of the resultant MILP problem. In fact, when the number of operations increases, the additional savings on generation cost usually become less significant; cf. Fig. \ref{fig:cost}. 
{\yuqi In addition, excessive switching operations can lead to security and stability issues in transmission grids  (see e.g.,\cite{fisher2008optimal}). Therefore, in practice the number of operations is normally restricted to a small number.} 
%{\hao (the former sentence reads really bad.)} 
In the tested 118-bus system, until $s=5$ operations, the proposed control scheme only requires less than $10$ seconds to find the solution.
The results imply the efficiency and scalability of the proposed optimization formulation for real-time implementation.

\subsubsection*{AC feasibility} It is certainly important to ensure that the {\yuqi switching decisions} of the topology optimization problem \eqref{eq:TS2} are feasible considering the ac power model. Thus, we verified the decisions of the breaker-level topology optimization shown in Table \ref{tab:sd} using the ac power flow. Our results confirm that the ac power flow model remains feasible in the post-event systems for the proposed topology optimization.

\section{CONCLUSIONS} \label{sec:con}

In this paper, we present a post-event analysis of substation bus splits and arrive at an equivalent bus-branch model for such events. 
Using this equivalent model, we propose a substation-level network topology optimization formulation that can incorporate both line switching and bus splitting. To deal with the bilinearity in the formulation, the McCormick relaxation has been utilized to devise a tractable MILP reformulation, which can be efficiently solved for real-time applications. Numerical studies on the IEEE 118-bus system corroborate the efficacy of the proposed topology optimization algorithm in terms of operational cost reduction and computational complexity. 
{\yuqi Future work includes topology optimization under an increased set of bus split operations such as the reconnection of multiple lines/injections. To accelerate the running time of ac power flow based topology optimization, we are also interested to explore a learning-based framework to approach the associated mixed-integer program.}

% \addtolength{\textheight}{-12cm}   % This command serves to balance the column lengths
                                  % on the last page of the document manually. It shortens
                                  % the textheight of the last page by a suitable amount.
                                  % This command does not take effect until the next page
                                  % so it should come on the page before the last. Make
                                  % sure that you do not shorten the textheight too much.

%%%%%%%%%%%%%%%%%%%%%%%%%%%%%%%%%%%%%%%%%%%%%%%%%%%%%%%%%%%%%%%%%%%%%%%%%%%%%%%%

%%%%%%%%%%%%%%%%%%%%%%%%%%%%%%%%%%%%%%%%%%%%%%%%%%%%%%%%%%%%%%%%%%%%%%%%%%%%%%%%

\bibliography{bibliography.bib}

% Generated by IEEEtran.bst, version: 1.14 (2015/08/26)
\begin{thebibliography}{10}
\providecommand{\url}[1]{#1}
\csname url@samestyle\endcsname
\providecommand{\newblock}{\relax}
\providecommand{\bibinfo}[2]{#2}
\providecommand{\BIBentrySTDinterwordspacing}{\spaceskip=0pt\relax}
\providecommand{\BIBentryALTinterwordstretchfactor}{4}
\providecommand{\BIBentryALTinterwordspacing}{\spaceskip=\fontdimen2\font plus
\BIBentryALTinterwordstretchfactor\fontdimen3\font minus
  \fontdimen4\font\relax}
\providecommand{\BIBforeignlanguage}[2]{{%
\expandafter\ifx\csname l@#1\endcsname\relax
\typeout{** WARNING: IEEEtran.bst: No hyphenation pattern has been}%
\typeout{** loaded for the language `#1'. Using the pattern for}%
\typeout{** the default language instead.}%
\else
\language=\csname l@#1\endcsname
\fi
#2}}
\providecommand{\BIBdecl}{\relax}
\BIBdecl

\bibitem{wood2013power}
A.~J. Wood, B.~F. Wollenberg, and G.~B. Shebl{\'e}, \emph{Power generation,
  operation, and control}.\hskip 1em plus 0.5em minus 0.4em\relax John Wiley \&
  Sons, 2013.

\bibitem{fisher2008optimal}
E.~B. Fisher, R.~P. O'Neill, and M.~C. Ferris, ``Optimal transmission
  switching,'' \emph{IEEE Trans. Power Systems}, vol.~23, no.~3, pp.
  1346--1355, 2008.

\bibitem{hedman2009optimal}
K.~W. Hedman, R.~P. O'Neill, E.~B. Fisher, and S.~S. Oren, ``Optimal
  transmission switching with contingency analysis,'' \emph{IEEE Trans. Power
  Systems}, vol.~24, no.~3, pp. 1577--1586, 2009.

\bibitem{ruiz2012tractable}
P.~A. Ruiz, J.~M. Foster, A.~Rudkevich, and M.~C. Caramanis, ``Tractable
  transmission topology control using sensitivity analysis,'' \emph{IEEE Trans.
  Power Systems}, vol.~27, no.~3, pp. 1550--1559, 2012.

\bibitem{qiu2015chance}
F.~Qiu and J.~Wang, ``Chance-constrained transmission switching with guaranteed
  wind power utilization,'' \emph{IEEE Trans. Power Systems}, vol.~30, no.~3,
  pp. 1270--1278, 2015.

\bibitem{zhou2020transmission}
Y.~Zhou, H.~Zhu, and G.~A. Hanasusanto, ``Transmission switching under wind
  uncertainty using linear decision rules,'' in \emph{Proc. IEEE PES General
  Meeting}, 2020.

\bibitem{mazi1986corrective}
A.~A. Mazi, B.~F. Wollenberg, and M.~H. Hesse, ``Corrective control of power
  system flows by line and bus-bar switching,'' \emph{IEEE Trans. Power
  Systems}, vol.~1, no.~3, pp. 258--264, 1986.

\bibitem{shao2005corrective}
W.~Shao and V.~Vittal, ``Corrective switching algorithm for relieving overloads
  and voltage violations,'' \emph{IEEE Trans. Power Systems}, vol.~20, no.~4,
  pp. 1877--1885, 2005.

\bibitem{zaoui2005coupling}
F.~Zaoui, S.~Fliscounakis, and R.~Gonzalez, ``Coupling {OPF} and topology
  optimization for security purposes,'' in \emph{15th Power Systems Computation
  Conference}, 2005, pp. 22--26.

\bibitem{heidarifar2015network}
M.~Heidarifar and H.~Ghasemi, ``A network topology optimization model based on
  substation and node-breaker modeling,'' \emph{IEEE Trans. Power Systems},
  vol.~31, no.~1, pp. 247--255, 2015.

\bibitem{pradeep2011cim}
Y.~Pradeep, P.~Seshuraju, S.~A. Khaparde, and R.~K. Joshi, ``{CIM}-based
  connectivity model for bus-branch topology extraction and exchange,''
  \emph{IEEE Trans. Smart Grid}, vol.~2, no.~2, pp. 244--253, 2011.

\bibitem{park2019sparse}
B.~Park, J.~Holzer, and C.~L. DeMarco, ``A sparse tableau formulation for
  node-breaker representations in security-constrained optimal power flow,''
  \emph{IEEE Trans. Power Systems}, vol.~34, no.~1, pp. 637--647, 2019.

\bibitem{park2020optimal}
B.~Park and C.~L. Demarco, ``Optimal network topology for node-breaker
  representations with {AC} power flow constraints,'' \emph{IEEE Access},
  vol.~8, pp. 64\,347--64\,355, 2020.

\bibitem{mccormick1976computability}
G.~P. McCormick, ``Computability of global solutions to factorable nonconvex
  programs: {Part I-Convex underestimating problems},'' \emph{Mathematical
  programming}, vol.~10, no.~1, pp. 147--175, 1976.

\bibitem{kekatos2012joint}
V.~Kekatos and G.~B. Giannakis, ``Joint power system state estimation and
  breaker status identification,'' in \emph{Proc. North American Power Symp.},
  2012.

\bibitem{korres2006substation}
G.~Korres, P.~Katsikas, and G.~Chatzarakis, ``Substation topology
  identification in generalized state estimation,'' \emph{International Journal
  of Electrical Power \& Energy Systems}, vol.~28, no.~3, pp. 195--206, 2006.

\bibitem{deka2015one}
D.~Deka, R.~Baldick, and S.~Vishwanath, ``One breaker is enough: Hidden
  topology attacks on power grids,'' in \emph{Proc. IEEE PES General Meeting},
  2015.

\bibitem{ten2017impact}
C.-W. Ten, K.~Yamashita, Z.~Yang, A.~V. Vasilakos, and A.~Ginter, ``Impact
  assessment of hypothesized cyberattacks on interconnected bulk power
  systems,'' \emph{IEEE Trans. Smart Grid}, vol.~9, no.~5, 2017.

\bibitem{zhou2018false}
Y.~Zhou, J.~Cisneros-Saldana, and L.~Xie, ``False analog data injection attack
  towards topology errors: Formulation and feasibility analysis,'' in
  \emph{Proc. IEEE PES General Meeting}, 2018.

\bibitem{jahromi2019cyber}
A.~A. Jahromi, A.~Kemmeugne, D.~Kundur, and A.~Haddadi, ``Cyber-physical
  attacks targeting communication-assisted protection schemes,'' \emph{IEEE
  Trans. Power Systems}, vol.~35, no.~1, pp. 440--450, 2019.

\bibitem{zhou2019bus}
Y.~Zhou and H.~Zhu, ``Bus split sensitivity analysis for enhanced security in
  power system operations,'' in \emph{Proc. North American Power Symp.}, 2019.

\bibitem{soroush2013accuracies}
M.~Soroush and J.~D. Fuller, ``Accuracies of optimal transmission switching
  heuristics based on {DCOPF} and {ACOPF},'' \emph{IEEE Trans. Power Systems},
  vol.~29, no.~2, pp. 924--932, 2013.

\bibitem{bai2016two}
Y.~Bai, H.~Zhong, Q.~Xia, and C.~Kang, ``A two-level approach to {AC} optimal
  transmission switching with an accelerating technique,'' \emph{IEEE Trans.
  Power Systems}, vol.~32, no.~2, pp. 1616--1625, 2016.

\bibitem{zamzam2020learning}
A.~S. Zamzam and K.~Baker, ``Learning optimal solutions for extremely fast {AC}
  optimal power flow,'' in \emph{2020 IEEE International Conference on
  Communications, Control, and Computing Technologies for Smart Grids
  (SmartGridComm)}.\hskip 1em plus 0.5em minus 0.4em\relax IEEE, 2020, pp.
  1--6.

\bibitem{kocuk2017new}
B.~Kocuk, S.~S. Dey, and X.~A. Sun, ``New formulation and strong {MISOCP}
  relaxations for {AC} optimal transmission switching problem,'' \emph{IEEE
  Trans. Power Systems}, vol.~32, no.~6, pp. 4161--4170, 2017.

\bibitem{griva2009linear}
I.~Griva, S.~G. Nash, and A.~Sofer, \emph{Linear and nonlinear
  optimization}.\hskip 1em plus 0.5em minus 0.4em\relax Siam, 2009, vol. 108.

\bibitem{gupte2013solving}
A.~Gupte, S.~Ahmed, M.~S. Cheon, and S.~Dey, ``Solving mixed integer bilinear
  problems using {MILP} formulations,'' \emph{SIAM Journal on Optimization},
  vol.~23, no.~2, pp. 721--744, 2013.

\bibitem{bhela2019designing}
S.~Bhela, D.~Deka, H.~Nagarajan, and V.~Kekatos, ``Designing power grid
  topologies for minimizing network disturbances: An exact {MILP}
  formulation,'' in \emph{Proc. American Control Conference (ACC)}, 2019.

\bibitem{bazrafshan2019optimal}
M.~Bazrafshan, N.~Gatsis, and H.~Zhu, ``Optimal power flow with step-voltage
  regulators in multi-phase distribution networks,'' \emph{IEEE Trans. Power
  Systems}, vol.~34, no.~6, pp. 4228--4239, 2019.

\end{thebibliography}

%\begin{comment}
%\begin{thebibliography}{1}

%\itemsep 2pt

%\bibitem{ANSI}
%emph{American National Standard
%For Electric Power Systems and Equipment- Voltage Ratings (60 Hertz)}, ANSI C84.1-2011, American National Standards Institute (ANSI), Inc. 

%\end{thebibliography}
%\end{comment}

\bibliographystyle{IEEEtran}

\end{document}